\documentclass[11pt, reqno]{amsart}

\title[Sign changes]{On signs of Fourier coefficients on $\GL(n)$}
\author{Didier Lesesvre, Ming Ho Ng, Yingnan Wang}
\date{\today}

\address{Universit\'e de Lille -- Laboratoire Paul Painlev\'e, UMR 8524,
 59000 Lille, France}
 \email{didier.lesesvre@univ-lille.fr}
 
 \address{Department of Mathematics, The Chinese University Of Hong Kong, Shatin, Hong Kong, P.R. China}
 \email{mhng@math.cuhk.edu.hk} 
 
 \address{School of Mathematical Sciences, Shenzhen University,
Shenzhen, Guangdong 518060, P.R. China}
 \email{ynwang@szu.edu.cn}

\usepackage{mathrsfs}
\usepackage{amssymb,amsmath,amsthm,color}
\usepackage[colorlinks=true,
linkcolor=blue, citecolor=blue]{hyperref}

\newtheorem{theorem}{Theorem}[section]
\newtheorem{lemma}{Lemma}[section]

\allowdisplaybreaks

\theoremstyle{definition}

\theoremstyle{remark}

\newcommand{\eval}[2][\right]{\relax
  \ifx#1\right\relax \left.\fi#2#1\rvert}

\newcommand{\GL}{\mathrm{GL}}
\newcommand{\SL}{\mathrm{SL}}
\renewcommand{\leq}{\leqslant}
\renewcommand{\geq}{\geqslant}
\renewcommand{\le}{\leqslant}
\renewcommand{\ge}{\geqslant}

\numberwithin{equation}{section}

\begin{document}

\maketitle

\begin{abstract}
We study statistical properties of Fourier coefficients of automorphic forms on $\GL(n)$. For most Hecke-Maass cusp forms, we give the asymptotic number of nonvanishing coefficients, show that there is a positive proportion of sign changes among them, when these are real, and describe the asymptotic density of these signs. We generalize the results by Jääsaari obtained in the case of self-dual forms of $\GL(3)$ and our method moreover circumvents the assumption of the Generalized Ramanujan Conjecture.
\end{abstract}

\section{Background}

Automorphic forms arise in various fashion, ranging from arithmetic to spectral theory, from geometry to representation theory. Langlands conjectures postulate they all stem from automorphic forms on the general linear groups $\mathrm{GL}(n)$, giving a particular importance to the understanding of these. Information about automorphic forms is often obtained by studying their Fourier coefficients, and we investigate in this paper statistical properties on the vanishing and the signs of the Fourier coefficients of automorphic forms on~$\mathrm{GL}(n)$.

Let $\mathcal{H}$ be an orthonormal basis of Hecke-Maass cusp forms \cite{goldfeld} for
$\SL(n,\mathbb{Z})$ with $n\ge 3$. For $T>1$, define
$$
\mathcal{H}_T=\left\{\phi\in\mathcal{H}:\mu_\phi\in i\mathbb{R}^n,\,\|\mu_\phi\|_2\leq T\right\}
$$
where $\mu_\phi \in \mathbb{C}^n$ is the Langlands parameters of $\phi \in\mathcal{H}$ and $\|\cdot\|_2$ is the Euclidean norm. The Weyl law \cite{mueller} for $\SL(n)$ states that $\#\mathcal{H}_T\asymp T^d$, where $d=n(n+1)/2-1$. 
A fundamental invariant attached to $\phi \in \mathcal{H}$ is its Godement-Jacquet $L$-function \cite{goldfeld}, which can be written
\begin{equation}
L(s, \phi) = \sum_{m=1}^\infty \frac{A_\phi(m, 1, \ldots, 1)}{m^s} = \prod_p \prod_{j=1}^n (1-\pi_{\phi, j}(p)p^{-s})^{-1}, 
\end{equation}
for $\Re(s) \gg 1$, where the $A_\phi(m_1, \ldots, m_{n-1})$ are called the Fourier coefficients of $\phi$, and the $\pi_{\phi, j}(p) \in \mathbb{C}$ are its Satake parameters at $p$. It is worth to note that the Fourier coefficient $A_\phi(m_1, \ldots, m_{n-1})$ is multiplicative. The fact that the forms have the trivial central character rephrases as $\pi_{\phi,1}(p)\pi_{\phi,2}(p)\cdots\pi_{\phi,n}(p)=1$, and the Generalized Ramanujan Conjecture (GRC) states that $|\pi_{\phi, j}(p)| = 1$ for all $1 \leqslant j \leqslant n$. 
This motivates to parametrize the Satake parameters by 
$$
\pi_{\phi,j}(p)=e^{i\theta_{\phi,j}(p)}
$$
where $\theta_{\phi,j}(p)\in [0,2\pi)\cup i\mathbb{R}\cup\pi +i\mathbb{R}$.

For $\mathbf{\kappa}=(\kappa_1,\ldots,\kappa_{n-1})\in \mathbb{N}_0^{n-1}$, the degenerate Schur polynomial is defined as
\begin{equation}\label{eq-schur}
S_{\mathbf{\kappa}} (x_{1}, x_{2},\ldots, x_{n}) := \frac{\det \begin{pmatrix} x_j^{\sum_{l=1}^{n-i} (\kappa_l +1)}\end{pmatrix}_{1\le i,j\le n}}{\det \begin{pmatrix} x_j^{\sum_{l=1}^{n-i} 1}\end{pmatrix}_{1\le i,j\le n}}
\end{equation}
and we have 
\begin{equation}\label{fourier-schur}
A_\phi(p^\kappa) := A_\phi(p^{\kappa_1}, \ldots, p^{\kappa_{n-1}})=S_{\mathbf{\kappa}} (\pi_{\phi,1}(p), \pi_{\phi,2}(p),\cdots, \pi_{\phi,n}(p)).
\end{equation}
If $x_j\in \mathbb{C}$ is such that $|x_j|=1$ for all $1\leq j\leq n$, then there exist constants $a_\mathbf{\kappa}<b_\mathbf{\kappa}$ only depending on $\mathbf{\kappa}$ and $n$ such that
$$
a_\mathbf{\kappa}\leq |S_\mathbf{\kappa}(x_{1}, x_{2},\ldots, x_{n})| \leq b_\mathbf{\kappa}.
$$
When $A_\phi(m_{n-1}, \ldots, m_1) = \overline{A_\phi(m_1, \ldots, m_{n-1})}$, we obtain that $A_\phi(p^{\kappa}) \in \mathbb{R}$ if we have $(\kappa_1,\ldots,\kappa_{n-1})=(\kappa_{n-1},\ldots,\kappa_{1})$.
%We assume this condition from now on, to be able to study the signs of $A_\phi(m_{n-1}, \ldots, m_1)$. 
Moreover, if the Generalized Ramanujan Conjecture holds at $p$, then for $(\kappa_1,\ldots,\kappa_{n-1})=(\kappa_{n-1},\ldots,\kappa_{1})$, we have
$$
a_\mathbf{\kappa}\leq A_\phi(p^{\kappa})\leq b_\mathbf{\kappa}.
$$

The Sato-Tate conjecture postulates the asymptotic distribution of the Satake parameters for each form $\phi \in \mathcal{H}$ when $p$ grows. Even though this horizontal statement is out of reach for general automorphic forms, a vertical version of the Sato-Tate conjecture --- where $p$ remains fixed but the form varies --- was established by Matz and Templier \cite{MT}. Further,  \cite{LNW} obtained a quantitative version with an explicit rate of convergence, which states that 
\begin{equation}
\frac{1}{|\mathcal{H}_T|} \left| \left\{ \phi \in \mathcal{H}_T \ : \ \underline{\theta}_\phi(p) \in I/\mathfrak{S}_n \right\} \right| = \int_{I/\mathfrak{S}_n} d\mu_p + O\left( \frac{\log p}{\log T}\right), 
\end{equation}
where $I = \prod_j [a_j, b_j]$ is a box in $[0, 2\pi)^n$ and the $p$-adic Plancherel measure is defined by
\begin{equation}
d\mu_p(\underline{\theta}) = \prod_{j=2}^n \frac{1-p^{-j}}{1-p^{-1}} \prod_{1 \leqslant \ell < m \leqslant n} |e^{i\theta_\ell} - p^{-1} e^{i\theta_m}|^{-2} d\mu_{\rm ST}(\underline{\theta})
\end{equation}
where the Sato-Tate measure is given by
\begin{equation}
d\mu_{\rm ST}(\underline{\theta}) = \frac{1}{n! (2\pi)^{n-1}} \prod_{1 \leqslant \ell < m \leqslant n} |e^{i\theta_\ell} - e^{i\theta_m}|^2 d\theta_1 \cdots d\theta_n.
\end{equation}

\section{Number of nonzero coefficients}

The following theorem is the analogue of \cite[Theorem 4]{Ja}; it gives the asymptotic number of nonvanishing coefficients $A_\phi(m^{\kappa})$ for a Hecke-Maass cusp form $\phi$. Here we do not require $A_\phi(m^{\kappa})$ to be real and they are arbitrary complex numbers. Moreover, by using knowledge about the exceptional set of automorphic forms, i.e. those violating the Generalized Ramanujan Conjecture, we are able to state an unconditional result.

\begin{theorem}\label{thm 2.1}
Let $c\log T\leq X\leq e^{(\log T)/k(T)}$ for a constant $c>0$, where $k(T)=o(\log T)$ and $k(T)/\sqrt{\log T} \to \infty$ as $T\rightarrow\infty$.
%Fix $\mathbf{\kappa}=(\kappa_1,\ldots,\kappa_{n-1})=(\kappa_{n-1},\ldots,\kappa_{1})$ in~$\mathbb{N}^{n-1}$.
As $T\rightarrow\infty$, there exists $\mathcal{H}_{T,\mathbf{\kappa}}\subset\mathcal{H}_T$ 
such that for all Hecke-Maass cusp forms $\phi\in\mathcal{H}_{T,\mathbf{\kappa}}$, we have
$$
\#\left\{m\leq X: A_\phi(m^\kappa) = A_\phi(m^{\kappa_1}, \ldots, m^{\kappa_{n-1}})\neq0\right\}\asymp X\prod_{p\leq X\atop A_\phi(p^\kappa)=0}\bigg(1-\frac1p\bigg)
$$
and, letting $d = \tfrac{1}{2}n(n+1)-1$, 
$$
\#(\mathcal{H}_{T}\backslash\mathcal{H}_{T,\mathbf{\kappa}})\ll T^de^{-(c_{\mathbf{\kappa},1}\log T)/\log X}+T^de^{-c_{\mathbf{\kappa},2}\log X}
$$ for some constants $c_{\mathbf{\kappa},1}$ and $c_{\mathbf{\kappa},2}$.
\end{theorem}

\begin{proof}
The upper bound is a classical result valid in a general setting under mild assumptions, as a consequence of sieve methods (see \cite{GKM}); we briefly recall the argument for completeness. Put 
$$\mathcal{B} = \{p \in \mathbb{P} \ : \ A_\phi(p^\kappa) = 0\}.$$
The set of integers $m \leqslant X$ such that $A_\phi(m^\kappa) \neq 0$ is covered by two sets:
\begin{align}
    S_1 & = \left\{m \leqslant X \ : \ p < X \text{ and } p\in \mathcal{B}\implies \ p \nmid m\right\}, \\
    S_2 & =  \left\{m\leqslant X \ : \ (m,P_X)>1\ \textrm{and}\ (m,P_X)^2\mid m\right\}, 
\end{align}
where we let 
\begin{equation}
P_X=\prod_{p\leq X\atop p\in \mathcal{B}}p.
\end{equation}
Indeed, by multiplicativity of Fourier coefficients and by factoring $m$ into product of primes, $A_\phi(m^\kappa) \neq 0$ only if $m$ contains no primes in $\mathcal{B}$ (which corresponds to $m \in S_1$), or for any prime $p$ in $\mathcal{B}$ which divides $m$, $p^2\mid m$,   (which corresponds to $m \in S_2$).

We use the sieve method to control the size of both sets, starting with $S_1$.
Introduce $\mathcal{A} = \{m \in \mathbb{N} \ : \ m \leqslant X\}$ and define the sifting set
\begin{equation}
S(\mathcal{A}, \mathcal{B}, z) := \# \{m \in \mathcal{A} \ : \ \forall p \in \mathcal{B}, \ p < z \implies \ p \nmid m\}, 
\end{equation}
for a parameter $z\geqslant 1$. This set is exactly $S_1$ when $z=X$.
We explain how the claimed upper bound is a consequence of \cite[Theorem 2.2]{Hal}, of which we follow the notations. Set $\mathcal{A}_d := \{m \in \mathcal{A} \ : \ m \equiv 0 \mod d\}$ for integers $d \geqslant 1$, and introduce
\begin{equation}
R_d := |\mathcal{A}_d| - \frac{1}{d} X.
\end{equation}
Let $\omega(p) = \mathbf{1}_{p \in \mathcal{B}}$ and $\omega(d) := \prod_{p\mid d} \omega(p)$, we verify the three hypotheses (denoted by $(\Omega_1)$, $(\Omega_2(\kappa))$, $(R)$ therein) required to apply~\cite[Theorem 2.2]{Hal}: first of all,
\begin{align}
\sum_{w < p < z} \frac{\omega(p) \log p}{p} \leqslant A\log \frac{z}{w} + A
\end{align}
for a certain $A>0$, as a direct consequence of Mertens' estimates, and this fulfills hypothesis $(\Omega_{2}(\kappa))$ in \cite{Hal}; moreover
\begin{equation}
0 \leqslant \frac{\omega(p)}{p} \leqslant 1-\frac{1}{A}
\end{equation}
for any $A>2$; which implies hypothesis $(\Omega_1)$ in \cite{Hal}; and finally
\begin{equation}
|R_d| \leqslant \omega(d)
\end{equation}
for all squarefree $d$ with no prime factor in $\overline{\mathcal{B}}$ (this is exactly \cite[Example 1 with $y=X$]{Hal}), and this is the content of hypothesis $(R)$ in \cite{Hal}. Therefore  \cite[Theorem 2.2]{Hal} ensures that 
\begin{equation}\label{set contains no p in B}
S(\mathcal{A}, \mathcal{B}, z) \ll X\prod_{p<z} \left( 1-\frac{\omega(p)}{p}\right).
\end{equation}
for all $z \leqslant X$, which in particular is of the claimed order of magnitude as in Theorem~\ref{thm 2.1} by taking $z=X$.

It remains to take into account the elements in $S_2$, i.e. including those $m \leqslant X$ such that $A_\phi(m^\kappa)\neq0$ which are not grasped in $S(\mathcal{A}, \mathcal{B}, X)$. These correspond to positive integers~$m$ which contain high powers of primes in $\mathcal{B}$.
Applying \cite[Theorem 2.2]{Hal} again, the carninality of $S_2$ is bounded by (letting $d=(m,P_X)$)
\begin{align}\label{set contains p in B}
\sum_{d\mid P_X}\sum_{0<m\leq X/d^2\atop (m,P_X/d)=1}1
&\ll \sum_{d\mid P_X}\frac{X}{d^2}\prod_{p\mid  P_X/d}\bigg(1-\frac1p\bigg)\nonumber\\
&\ll X\prod_{p\mid P_X}\bigg(1-\frac1p\bigg)\sum_{d\mid P_X}\prod_{p\mid d}\frac1{p^2-p}\nonumber\\
&\ll X\prod_{p<X\atop p\in\mathcal{B}} \left( 1-\frac{1}{p}\right).
\end{align}

Combining \eqref{set contains no p in B} and \eqref{set contains p in B}, we obtain the upper bound claimed in Theorem \ref{thm 2.1}.

We turn to the proof of the lower bound. We only prove the case that $\mathbf{\kappa}=(\kappa_1,\ldots,\kappa_{n-1})\neq(\kappa_{n-1},\ldots,\kappa_{1})=:\mathbf{\kappa}^\iota$. Otherwise, the proof is similar with obvious modifications. 

Define $\|\mathbf{\kappa}\|:=\sum_{j=1}^{n-1}(n-j)\kappa_j$ and $|\mathbf{\kappa}|=\sum_{j=1}^{n-1}\kappa_j$. If $\mathbf{\kappa}\neq\mathbf{\kappa}^\iota$, by \cite[(2.3)]{LNRW}, we have
\begin{equation}\label{Hecke}
A_\phi(p^\kappa)^2=\sum\limits_{\mathbf{\xi}\neq\mathbf{0}\atop\|\mathbf{\xi}\|\le 2\|\mathbf{\kappa}\|}
d_{\mathbf{\kappa}\mathbf{\kappa}}^{\mathbf{\xi}}A_\phi(p^\xi).
\end{equation}
Here $\mathbf{\xi}\neq\mathbf{0}$ is assured by the fact that $\{S_\mathbf{\kappa}(e^{i\theta_1}, e^{i\theta_2},\ldots, e^{i\theta_n}):\mathbf{\kappa}=(\kappa_1,\ldots,\kappa_{n-1})\in \mathbb{N}_0^{n-1}\}$ forms an orthonormal set under the inner product with repect to the measure $d\mu_{\rm ST}(\underline{\theta})$, see \cite[(2.4)]{LNRW}.

Introduce the exceptional set
$$
\mathcal{E}(T,p)=\left\{\phi\in\mathcal{H}_T:\max_{1\leq \ell\leq n}|\pi_{\phi,\ell}(p)|>2\right\}.
$$
Then by \cite[Corollary 1.8]{MT}, $\#\mathcal{E}(T,p)\ll T^{d-c_0/\log p}$ for $p\leq X$, where $c_0$ is a constant depending only on $n$. Hence, there exits a constant $c_1$ only depending on $n$ such that
\begin{equation*}\label{exceptional set}
\bigcup_{p\leq X}\#\mathcal{E}(T,p)\ll e^{(\log T)/k(T)}T^{d-c_0/\log X}\ll T^d e^{-c_1(\log T)/\log X}.
\end{equation*}

For any $\phi\in\mathcal{H}_T\backslash\bigcup_{p\leq X}\#\mathcal{E}(T,p)$, we have $|A_\phi(p^\mathbf{\kappa})|\leq B_{\mathbf{\kappa}}$ for some positive constant $B_{\mathbf{\kappa}}$ only depending on $\mathbf{\kappa}$.
Then, by \eqref{Hecke} and the Hecke relations \cite[(4.5)]{LNRW} with $\mathbf{\kappa}=(\kappa_1,\ldots,\kappa_{n-1})$, we have 
\begin{align*}
& \sum_{p\leq X\atop A_\phi(p^\mathbf{\kappa})\neq0}\frac{1}{p} \geq \sum_{\log^{2025} X\leq p\leq X^{1/2025}\atop \mathrm{Re}A_\phi(p^\mathbf{\kappa})<0}\frac{1}{p} \nonumber\\
& \quad \geqslant \sum_{\log^{2025} X\leq p\leq X^{1/2025}}\frac{(\mathrm{Re} A_\phi(p^\kappa))^2-B_{\mathbf{\kappa}}\mathrm{Re}A_\phi(p^\kappa)}{2B_{\mathbf{\kappa}}^2p}\nonumber\\
& \quad \geq\sum_{\log^{2025} X\leq p\leq X^{1/2025}}\frac{ A_\phi(p^\kappa)^2 + \overline{A_\phi(p^\kappa)}^2 + 2|A_\phi(p^\kappa)|^2 -2B_{\mathbf{\kappa}}A_\phi(p^\kappa)- 2B_{\mathbf{\kappa}}\overline{A_\phi(p^\kappa)}}{8B_{\mathbf{\kappa}}^2p} \nonumber\\
& \quad =\sum_{\log^{2025} X\leq p\leq X^{1/2025}}\frac{\sum\limits_{\mathbf{\xi_1}\neq\mathbf{0}\atop\|\mathbf{\xi_1}\|\le 2\|\mathbf{\kappa}\|}
d_{\mathbf{\kappa}\mathbf{\kappa}}^{\mathbf{\xi_1}}A_\phi(p^{\mathbf{\xi_1}}) + \sum\limits_{\mathbf{\xi_2}\neq\mathbf{0}\atop\|\mathbf{\xi_2}\|\le 2\|\mathbf{\kappa}^\iota\|}
d_{\mathbf{\kappa}^\iota\mathbf{\kappa}^\iota}^{\mathbf{\xi_2}}A_\phi(p^{\mathbf{\xi_2}})}{8B_{\mathbf{\kappa}}^2p} \nonumber\\
& \quad \qquad +\sum_{\log^{2025} X\leq p\leq X^{1/2025}}\frac{ 2+2\sum\limits_{\mathbf{\xi_3}\neq\mathbf{0}\atop\|\mathbf{\xi_3}\|\le n|\mathbf{\kappa}|}
d_{\mathbf{\kappa}\mathbf{\kappa}^\iota}^{\mathbf{\xi_3}}A_\phi(p^{\mathbf{\xi_3}})-2B_{\mathbf{\kappa}}A_\phi(p^\kappa) - 2B_{\mathbf{\kappa}}\overline{A_\phi(p^\kappa)}}{8B_{\mathbf{\kappa}}^2p} \nonumber\\
& \quad =\frac{1+o(1)}{4B_{\mathbf{\kappa}}^2}\log\log X +\sum_{\log^{2025} X\leq p\leq X^{1/2025}} \frac{\sum\limits_{\mathbf{\xi_1}\neq\mathbf{0}\atop\|\mathbf{\xi_1}\|\le 2\|\mathbf{\kappa}\|}
d_{\mathbf{\kappa}\mathbf{\kappa}}^{\mathbf{\xi_1}}A_\phi(p^{\mathbf{\xi_1}}) + \sum\limits_{\mathbf{\xi_2}\neq\mathbf{0}\atop\|\mathbf{\xi_2}\|\le 2\|\mathbf{\kappa}^\iota\|}
d_{\mathbf{\kappa}^\iota\mathbf{\kappa}^\iota}^{\mathbf{\xi_2}}A_\phi(p^{\mathbf{\xi_2}})}{8B_{\mathbf{\kappa}}^2p}\nonumber \nonumber\\
& \quad \qquad +\sum_{\log^{2025} X\leq p\leq X^{1/2025}} \frac{ \sum\limits_{\mathbf{\xi_3}\neq\mathbf{0}\atop\|\mathbf{\xi_3}\|\le n|\mathbf{\kappa}|}
d_{\mathbf{\kappa}\mathbf{\kappa}^\iota}^{\mathbf{\xi_3}}A_\phi(p^{\mathbf{\xi_3}})-2B_{\mathbf{\kappa}}A_\phi(p^\kappa) - 2B_{\mathbf{\kappa}}\overline{A_\phi(p^\kappa)}}{8B_{\mathbf{\kappa}}^2p}.
\end{align*}

%Hence
%\begin{align*}
%&\sum_{\phi\in \mathcal{H}_T\backslash\bigcup_{p\leq X}\#\mathcal{E}(T,p)}\sum_{p\leq X\atop A_\phi(p,1,\ldots,1,p)<0}\frac{1}{p}\nonumber\\
%&\geq \frac{1+o(1)}{4n^2}(\log\log X)\# \mathcal{H}_T\backslash\bigcup_{p\leq X}\#\mathcal{E}(T,p)\nonumber\\
%&\ +\sum_{\phi\in \mathcal{H}_T\backslash\bigcup_{p\leq X}\#\mathcal{E}(T,p)}\sum_{\log^{2025} X\leq p\leq X^{1/2025}}\frac{\sum\limits_{\mathbf{\xi}\neq\mathbf{0}\atop\|\mathbf{\xi}\|\le n|\mathbf{\kappa}|}
%d_{\mathbf{\kappa}\mathbf{\kappa}^\iota}^{\mathbf{\xi}}A_\phi(p^{\xi_1},\ldots,p^{\xi_{n-1}})-(4n^2-1)A_\phi(p,1,\ldots,1,p)}{4n^2p}.
%\end{align*}

On the other hand, by \cite[Theorem 1.1]{LNRW} with $j=[(\log X)/2025\log\log X]$, we have
\begin{align*}
&\sum_{\phi\in \mathcal{H}_T\backslash\bigcup_{p\leq X}\#\mathcal{E}(T,p)}\left|\sum_{\log^{2025} X\leq p\leq X^{1/2025}}\frac{
A_\phi(p^\xi)}{p}\right|^{2j}\\
&\ll (\log X)^{2j}\sum_{0\leq \ell\leq \log X}\sum_{\phi\in \mathcal{H}_T}\left|\sum_{2^\ell\log^{2025} X\leq p\leq 2^{\ell+1}\log X}\frac{
A_\phi(p^\xi)}{p}\right|^{2j}\\
&\ll (\log X)^{2j}\sum_{0\leq \ell\leq \log X}\left(T^d\left(\frac{C_{\mathbf{\kappa}}^2 j}{2^\ell\log^{2025} X\log (2^\ell\log^{2025} X)}\right)^j
+T^{d-1/2}\left(\frac{C_{\mathbf{\kappa}}Q^{L\|\mathbf{\kappa}\|}}{\log (2^\ell\log^{2025} X)}\right)^{2j}\right).
\end{align*}
Define
$$
E(p^\xi)=\left\{\phi\in \mathcal{H}_T\backslash\bigcup_{p\leq X}\#\mathcal{E}(T,p):\left|\sum_{\log^{2025} X\leq p\leq X^{1/2025}}\frac{
A_\phi(p^\xi)}{p}\right|\gg \log\log X\right\}.
$$
Then $\#E(p^\xi)$ is bounded by
\begin{align*}
&(\log X)^{2j}\sum_{0\leq \ell\leq \log X}\left(T^d\left(\frac{C_{\mathbf{\kappa}}^2 j}{2^\ell\log^{2025} X\log (2^\ell\log^{2025} X)}\right)^j
+T^{d-1/2}\left(\frac{C_{\mathbf{\kappa}}Q^{L\|\mathbf{\kappa}\|}}{\log (2^\ell\log^{2025} X)}\right)^{2j}\right)\\
&\ll T^d e^{-c_2\log X}.
\end{align*}
Therefore, for any $\phi$ in the set
\begin{align*}
\mathcal{H}_{T,\mathbf{\kappa}}=
\bigg(\mathcal{H}_T\backslash\bigcup_{p\leq X}\#\mathcal{E}(T,p)\bigg)
&\backslash\bigg(\mathop{\bigcup}\limits_{\mathbf{\xi_1}\neq\mathbf{0}\atop\|\mathbf{\xi_1}\|\le 2\|\mathbf{\kappa}\|}E(p^{\mathbf{\xi_1}})\mathop{\bigcup}\limits_{\mathbf{\xi_2}\neq\mathbf{0}\atop\|\mathbf{\xi_2}\|\le 2\|\mathbf{\kappa}^\iota\|}E(p^{\mathbf{\xi_1}})\mathop{\bigcup}\limits_{\mathbf{\xi_3}\neq\mathbf{0}\atop\|\mathbf{\xi_3}\|\le n|\mathbf{\kappa}|}E(p^{\mathbf{\xi_3}})\bigg)\\
&\backslash\bigg(\bigcup E(p^\mathbf{\kappa})\bigcup E(p^{\mathbf{\kappa}^\iota})\bigg),
\end{align*}
we have
\begin{align}\label{lower bound of negative primes}
\sum_{p\leq X\atop \mathrm{Re}A_\phi(p^\mathbf{\kappa})<0}\frac{1}{p}\gg \log\log X.
\end{align}
The lower bound follows by \cite[Theorem 1]{GKM} with $\mathcal{P}=\{p\leq x:\mathrm{Re}A_\phi(p^\kappa)\neq0\}$ and  $\mathcal{E}=\{p\leq x:\mathrm{Re}A_\phi(p^\kappa)=0\}$.

%
%\begin{equation}
%\sum_{\substack{p > x \\ A_\phi(p^{\kappa_1},\ldots,p^{\kappa_{n-1}}) = 0}} \frac{1}{p} \longrightarrow 0 \qquad \text{and} \quad \sum_{\substack{p \leqslant x \\ A_\phi(p^{\kappa_1},\ldots,p^{\kappa_{n-1}}) < 0}} \frac{1}{p} \longrightarrow \infty.
%\end{equation}
%This, in addition to the multiplicativity of the coefficients, allows to appeal to \cite[Lemma 2.4]{MR2} and conclude that 
%$$
%\#\left\{m\leq X: A_\phi(m,1,\ldots,1,m)\neq0\right\}\asymp X\prod_{p\leq X\atop A_\phi(p,1,\ldots,1,p)=0}\bigg(1-\frac1p\bigg)
%$$
%finishing the proof.
\end{proof}

\section{Positive proportion of sign changes}

The following theorem is the analogue of \cite[Theorem 3]{Ja}; it states that there is a positive proportion of sign changes among the real Fourier coefficients of almost all forms, in a quantitative way.

\begin{theorem}
\label{thm:proportion1}
Let $\varepsilon>0$ and $\mathbf{\kappa} \in  \mathbb{N}^{n-1}$ such that $(A_\phi(m^\kappa))_{m \geqslant 1}$ is real. Then there exists a subset $\mathcal{S}_T\subset\mathcal{H}_T$  with at least $(1-\varepsilon)\#\mathcal{H}_T$ elements such that for any $\phi\in\mathcal{S}_T$ the sequence $\{A_\phi(m^\kappa)\}_{m\geqslant 1}$ has a positive proportion of sign changes as $T\rightarrow\infty$.
\end{theorem}

\noindent \textit{Remark.} For instance, assuming $\mathbf{\kappa}=(\kappa_1,\ldots,\kappa_{n-1})=(\kappa_{n-1},\ldots,\kappa_{1})$ ensures that the coefficients $A_\phi(m^\kappa)$ are real.

We will appeal to the following lemma controlling the size of the pre-image for which corresponds to the small values of a mulivariate polynomial on a product of unit circles. 

\begin{lemma}
\label{lem:polynomial-small-values}
Let $n \geqslant 1$.
Let $P \in \mathbb{C}[X_1, \ldots, X_n]$ be a non-constant complex polynomial in $n$ variables. For all $\delta > 0$, the measure of the set
\begin{equation}
\left\{ (\theta_1, \ldots, \theta_n) \in [0, 2\pi]^n \ : \ \left|P\left(e^{i\theta_1}, \ldots, e^{i\theta_n}\right)\right| < \delta  \right\}
\end{equation}
is bounded by $\delta^{1/2^n}$.
\end{lemma}

\begin{proof}
We proceed by induction on the number of variables $n$, the result being straightforward for $n=1$ where in this case the polynomials are linear.
Assume $n \geqslant 2$ and denote the algebraic form of $P$ by 
\begin{equation}
P = \sum_{t_1=0}^{\gamma_{1}}\cdots\sum_{t_n=0}^{\gamma_{n}}a_{t_1\ldots t_n} X_1^{t_1} \cdots X_n^{t_n}, 
\end{equation}
where the $a_{t_1\ldots t_n} \in \mathbb{C}$ are the coefficients of $P$ and $\gamma_1, \ldots, \gamma_n \geqslant 0$ are the degrees associated with $X_1, \ldots, X_n$ respectively. Since $P$ is non-constant, we can assume $\gamma_n \geqslant 1$. We can therefore write
\begin{align*}
 P(e^{i\theta_1}, \ldots, e^{i\theta_n}) &= \sum_{t_1=0}^{\gamma_{1}}\cdots\sum_{t_n=0}^{\gamma_{n}}a_{t_1\ldots t_n}e^{it_1\theta_{1}}\cdots e^{it_n\theta_{n}} =\sum_{k=0}^{\gamma_{n}}b_k(e^{i\theta_1},\ldots,e^{i\theta_{n-1}}) e^{ik\theta_{n}}
\end{align*}
where we group all the terms indexing $t_1, \ldots, t_{n-1}$ into $b_{n}$, which is a polynomial of $n-1$ variables in $e^{i\theta_1}, \ldots, e^{i\theta_{n-1}}$. We decompose the interval $[0, 2\pi]$ in even subintervals of length $\delta$ of the form $[s_{\ell}\pi {\delta},(s_{\ell}+1)\pi {\delta}]$ for  $1\leq\ell\leq n-1$ and $0\leq s_{\ell}\leq 2/\delta-1$. 
For  each $1\leq\ell\leq n-1$, let $0\leq s_{\ell}\leq 2/\delta-1$ such that  $\theta_{\ell}\in[s_{\ell}\pi \delta,(s_{\ell}+1)\pi\delta ]$.  For $0\leq k\leq \gamma_{n}$, Taylor approximation ensures that
$$
b_{k}(e^{i\theta_{1}},\ldots,e^{i\theta_{n-1}})
=b_k(e^{is_{1}\pi \delta},\ldots,e^{is_{n-1}\pi \delta})+O({\delta}),
$$
where the implied constant only depends on $P$.
Therefore, $|P(e^{i\theta_1}, \ldots, e^{i\theta_n})|<\delta$ implies
$$
\left|\sum_{k=0}^{\gamma_{n}}b_k(e^{is_{1}\pi {\delta}},\ldots,e^{is_{n-1}\pi {\delta}}) e^{ik\theta_{n}}\right|\ll {\delta}, 
$$
i.e.
\begin{align*}
\left|b_{\gamma_{n}}(e^{ is_{1}\pi \delta},\ldots,e^{ is_{n-1}\pi \delta})\prod_{j=1}^{\gamma_{n}}(e^{i\theta_{n}}-z_j)\right|
\ll \delta,
\end{align*}
where the $z_j$ denote the roots of the polynomial $\sum_{k=0}^{\gamma_{n}}b_k(e^{is_{1}\pi \delta},\ldots,e^{is_{n-1}\pi \delta})z^k$ in $\mathbb{C}$.

If $b_{\gamma_{n}}(e^{ is_{1}\pi \delta},\ldots,e^{ is_{n-1}\pi \delta})\gg \delta^{1/2}$, then the product of the $e^{i\theta_{n}} - z_j$ is small and therefore $|e^{i\theta_{n}} - z_j| \ll \delta^{1/2}$ for at least one $j$. This implies that the values of $\theta_n$ are constrained in a set of length bounded by $\delta^{1/2}$ and the result follows. 
If $b_{\gamma_{n}}(e^{ is_{1}\pi \delta},\ldots,e^{ is_{n-1}\pi \delta})\ll \delta^{1/2}$, the result follows by induction since $b_{\gamma_n}$ is a polynomial in the  $n-1$ variables $e^{is_1}, \ldots, e^{is_{n-1}}$.
\end{proof}

\begin{proof}[Proof of Theorem \ref{thm:proportion1}]
By \cite[Corollary 3]{MR}, the sequence $\{A_\phi(m^\kappa)\}_{m\geqslant 1}$ has a positive proportion of sign changes as $T\rightarrow\infty$ if and only if $A_\phi(m^\kappa)<0$ for some $m \geqslant 1$ and $A_\phi(m^\kappa) \neq 0$ for a positive proportion of integers. 

By \cite[Theorem 1.5]{LNRW} with $\mathcal{P}=\{\textrm{all the primes}\}$, $\kappa=\{\kappa_1,\ldots,\kappa_{n-1}\}$ and $\delta=1$, there exits a subset $\mathcal{H}_{T}'\subset\mathcal{H}_T$ satisfying that for any $\phi\in\mathcal{H}_{T}'$, there exists a natural number~$m\geqslant 1$ such that $A_\phi(m^\kappa)<0$ and $\#(\mathcal{H}_T\backslash\mathcal{H}_{T}')\ll T^de^{-C_2\log T/\log_2T}$, where~$C_2$ and the implied constants depend on $\mathbf{\kappa}$ and $n$.

By \cite[Theorem 14]{Se}, to prove that $A_\phi(m^\kappa) \neq 0$ for a positive proportion of integers, it suffices to prove that
$$
\sum_{p\atop A_\phi(p^\kappa)=0}\frac1p<\infty.
$$

The vanishing of $A_\phi(p^\kappa)$ implies strong constraints on the parameters $\underline{\theta}_\phi$, and this cannot happen too often by the known average versions towards the Ramanujan conjecture~\cite{LNW}. More precisely, we have that for $\delta>0$, the set $$\left\{\phi\in\mathcal{H}_T:|A_\phi(p^\kappa)| = 0\right\}$$ is a subset of the union of \begin{equation}\label{bad primes}
\left\{\phi\in\mathcal{H}_T:\max\limits_{1\leq\ell\leq n}\log|\pi_{\phi,\ell}(p)|>0\right\}
\end{equation}
and
\begin{equation}\label{good primes}
\left\{\phi\in\mathcal{H}_T: \phi\text{ satisfies GRC and } |A_\phi(p^\kappa)|<p^{-\delta}\right\}.
\end{equation}
In the second case, since $A_\phi(p^\kappa)=S_\mathbf{\kappa}(e^{i\theta_{\phi,1}(p)},\ldots,e^{i\theta_{\phi,n}(p)})$ is a polynomial in the variables $(e^{i\theta_{\phi,1}(p)},\ldots,e^{i\theta_{\phi,n}(p)})$, Lemma \ref{lem:polynomial-small-values} ensures that the measure of the corresponding $\underline{\theta}_\phi \in [0,2\pi]^n$ is smaller than $p^{-\delta/2^n}$ for all $\delta > 0$.
Therefore, applying \cite[Theorem 1.1]{LNW} to~\eqref{bad primes} and applying \cite[Theorem 1.2, (5), (6)]{LNW} to \eqref{good primes}, we have that
$$
\frac{\#\left\{\phi\in\mathcal{H}_T:|A_\phi(p^\kappa)| = 0\right\}}{\#\mathcal{H}_T}\ll \frac1{p^{\delta/2^n}}+\frac{\log p}{\log T}.
$$
We deduce
$$
\sum_{\phi\in\mathcal{H}_T}\sum_{p\leq X\atop A_\phi(p^\kappa)=0}\frac1p\leq \#\mathcal{H}_T\sum_{p\leq X}\frac1p\bigg( \frac1{p^{\delta/2^n}}+\frac{\log p}{\log T}\bigg)\ll \#\mathcal{H}_T
$$
for any $X>1$ satisfying $\log X\ll\log T$ when $T$ is sufficiently large. By Chebyshev inequality, there hence exists an absolute constant $C>0$ and $\mathcal{S}_T'\subset \mathcal{H}_T$ such that, for all Hecke-Maass forms in $\mathcal{S}_T'$, we have
\begin{equation}\label{Ap=0}
\sum_{p\leq X\atop A_\phi(p^\kappa)=0}\frac1p\leq \frac{C}{2\varepsilon}
\end{equation}
and $\#(\mathcal{H}_T\backslash \mathcal{S}_T')\leq \varepsilon\#\mathcal{H}_T$. In particular, the series \eqref{Ap=0} converges as $T$ grows if we choose $X$ such that $X\rightarrow\infty$ as $T\rightarrow\infty$. By \cite[Theorem 14]{Se}, this implies that there is a positive proportion of $m$ such that $A_\phi(m^\kappa) \neq 0$, and hence concludes the proof by 
putting
$$
\mathcal{S}_T=\mathcal{S}_T'\cap \mathcal{H}_T'.
$$

\end{proof}

\section{Fourier coefficients of same signs}
The following theorem is the analogue of \cite[Theorem 6]{Ja}.
\begin{theorem}
\label{thm 4.1}
Let $\varepsilon>0$ and $\mathbf{\kappa} \in  \mathbb{N}^{n-1}$ such that $(A_\phi(m^\kappa))_{m \geqslant 1}$ is real. Then there exists a subset $\mathcal{S}_T''\subset\mathcal{H}_T$  with at least $(1-\varepsilon)\#\mathcal{H}_T$ elements such that for any $\phi\in\mathcal{S}_T''$ asymptotically half of non-zero coefficients $\{A_\phi(m^\kappa)\neq0\}_{m\geqslant 1}$ are positive and half of them are negative as $T\rightarrow\infty$.
\end{theorem}
\begin{proof}
Let $\mathcal{H}_{T,\mathbf{\kappa}}$ and $\mathcal{S}_T$ be defined as in Theorem \ref{thm 2.1} and Theorem \ref{thm:proportion1} respectively. Put
$$
\mathcal{S}_T''=\mathcal{H}_{T,\mathbf{\kappa}}\cap \mathcal{S}_T.
$$
Then $\mathcal{S}_T''$ contains at least $(1-\varepsilon)\#\mathcal{H}_T$ elements.
Furthermore, for any $\phi\in \mathcal{S}_T''$, by \eqref{lower bound of negative primes} and \eqref{Ap=0}, we have
$$
\sum_{p\leq X\atop A_\phi(p^\kappa)<0}\frac{1}{p}\gg \log\log X
$$
and
$$
\sum_{p\leq X\atop A_\phi(p^\kappa)=0}\frac1p\leq \frac{C}{2\varepsilon}.
$$
Then the theorem follows by \cite[Lemma 2.4]{MR2}.

\end{proof}

\subsection*{Acknowledgements} This work was done when D. L. was visiting Shenzhen University, which we thank for excellent working environment and support. D. L. acknowledges the support of the CDP C2EMPI, together with the French State under the France-2030 program, the University of Lille, the Initiative of Excellence of the University of Lille and the European Metropolis of Lille for their funding and support of the R-CDP-24-004-C2EMPI project. Y. W. is supported by National Natural Science Foundation of China (Grant No. 12371006). We moreover thank Jie Wu for enlightening discussions.

\end{document}